\documentclass{amsart}
\usepackage{amsthm}
\usepackage{tikz}
\usepackage{hyperref}
\usepackage{textcomp}
\usepackage{alltt}
\usepackage{changepage}
\usepackage[math]{cellspace}
\usepackage{url}
\usepackage[T1]{fontenc}
\usepackage[margin=1in]{geometry}
\theoremstyle{definition}
\newtheorem{conjecture}{Conjecture}

\newenvironment{absolutelynopagebreak}
{\par\nobreak\vfil\penalty0\vfilneg
\vtop\bgroup}
{\par\xdef\tpd{\the\prevdepth}\egroup
\prevdepth=\tpd}
\title{Betti numbers of unordered configuration spaces of small graphs}
\author{Gabriel C. Drummond-Cole}
\thanks{
I would like to thank Byung Hee An and Ben Knudsen for help and suggestions related to this project.
This work was supported by IBS-R003-D1. 
Last modified \today.
}
\begin{document}
\begin{abstract}
The purpose of this document is to provide data about known Betti numbers of unordered configuration spaces of small graphs in order to guide research and avoid duplicated effort.
It contains information for connected multigraphs having at most nine edges which contain no loops, no bivalent vertices, and no internal (i.e., non-leaf) bridges.
\end{abstract}
\maketitle
\section*{Introduction}
The purpose of this document is to provide data about known Betti numbers of unordered configuration spaces of small graphs in order to guide research and avoid duplicated effort.

It contains information for connected multigraphs having at most nine edges which contain no loops, no bivalent vertices, and no internal (i.e., non-leaf) bridges.
For each such graph there is a table of Betti numbers for configurations of small numbers of points along with Poincar\'e series which encode Betti numbers for all larger numbers of points. 
Because the Poincar\'e series are rational functions, this means that for fixed $i$ and $\Gamma$, the sequence $\beta_i (B_k\Gamma)$ of $i$th Betti numbers of is eventually polynomial in the number of points $k$ in the configuration. 
For the reader's convenience, the stable polynomial is also indicated.

The enumeration of the graphs in question was done using utilities that are part of Nauty~\cite{McKayPiperno:PCIII}.
For graphs with fewer than seven edges, this was also verified using Richard Mathar's preprint~\cite[Table 61]{Mathar:SSG}. 
The core homology computation was done using Macaulay2~\cite{GraysonStillman:M2}. 

The approach uses the fact that the chains of the unordered configuration spaces of a graph are a differential graded module over a polynomial ring with variables the edges of the graph.
The calculation is simplified further by a presentation of this differential graded module as the tensor product of finitely generated models for local configurations at the vertices of the graph.
This picture was pioneered by~\'Swi\k{a}tkowski~\cite{Swiatkowski:EHDCSG}; the action by the polynomial ring was described in~\cite{AnDrummond-ColeKnudsen:SSGBG} and analyzed in~\cite{AnDrummondColeKnudsen:ESHGBG}.
Other related work includes the following.
\begin{enumerate}
\item Ko and Park~\cite{KoPark:CGBG} calculated the first Betti numbers for all graphs.
\item Ramos~\cite{Ramos:SPHTBG} calculated all Betti numbers for trees. For trees all values are stable. In this document, this includes the graphs with fewer than two essential vertices.
\item Maci\k{a}\.{z}ek and Sawicki~\cite{MaciazekSawicki:NAQSG} calculated all Betti numbers for theta graphs. In this document, this data is reproduced in a different form in the section for graphs with precisely two essential vertices.
\end{enumerate}
A web version of this document is available \href{https://drummondcole.com/gabriel/academic/graph_configuration_betti_numbers/}{here}.

\subsection*{Examples of motivation}
Two examples of the use of this data are as follows.

In~\cite{AnDrummondColeKnudsen:ESHGBG}, we showed that the $i$th Betti numbers of $B_k(\Gamma)$ grows like a polynomial in $k$ of degree $\Delta_\Gamma^i-1$, where $\Delta_\Gamma^i$ is the maximal number of connected components of the complement of $i$ points in $\Gamma$. 
At the time, An made the following conjecture, which I believe has not appeared publicly, about the leading coefficient of this polynomial.
\begin{conjecture}[Leading coefficient conjecture]
For any graph $\Gamma$ and any index $i>1$, the leading term of the polynomial governing the growth of the $i$th Betti numbers of $B_k(\Gamma)$ is $C_\Gamma^i k^{\Delta_\Gamma^i-1}$, where
\[
C_\Gamma^i = \frac{1}{(\Delta_\Gamma^i-1)!}\sum_S \prod_{v\in S} (d(v)-2)
\]
where $S$ runs over $i$-element subsets of the essential vertices of $\Gamma$ and $d(v)$ is the degree of a vertex.
\end{conjecture}
This conjecture would be false for $i=1$ more or less because of the special case of Lemma 3.18 of~\cite{AnDrummondColeKnudsen:ESHGBG}. 

In any event, the data here provides evidence for the conjecture, which holds for all of these small graphs.

As a second example, as expressed above, for each $\Gamma$ and $i$ the sequence $\beta_i(B_k\Gamma)$ is eventually polynomial in $k$. 
A priori this does not tell us \emph{when} polynomiality starts. 
Consequently, an interesting question is whether it is possible to find a lower bound on the value of $k$ for which $\beta_i(B_k\Gamma)$ is polynomial in $k$ as a ``simple'' function of $i$ and ``simple'' invariants of $\Gamma$. 
For instance, is it possible that the beginning of the polynomial range is bounded below by a constant plus the number of essential vertices of $\Gamma$ for $\Gamma$ without isolated vertices?
The data here indicates that if so, the constant is at least two.

\subsection*{Reducing to the kinds of graphs considered here}
Betti numbers for unordered configuration spaces of graphs with loops can be obtained by reducing to loop-free graphs with the same number of edges~\cite[Lemma 4.6]{AnDrummondColeKnudsen:ESHGBG}. 
Betti numbers for unordered configuration spaces of loop-free graphs with bivalent vertices can be obtained from such information for graphs without bivalent vertices (smoothing bivalent vertices is a homeomorphism and thus does not affect the configuration spaces). 
Betti numbers for unordered configuration spaces of disconnected graphs can be obtained combinatorially from the counts for connected graphs by summing over all ways of partitioning the desired number of points among the path components (this is true for any locally path-connected topological spaces and is not particular to graphs). 
Betti numbers for unordered configuration spaces of graphs with internal bridges can be obtained in a similar combinatorial manner from such information for the graphs obtained by cutting the edge into two leaves~\cite[Proposition 5.22]{AnDrummond-ColeKnudsen:SSGBG}

\subsection*{Notes about the data}
The graphs are organized by the number of essential (valence at least three) vertices. 
They are presented in sparse6 format, with degree sequences, with adjacency matrices, and with a visual representation.
The number of essential vertices also gives the maximal homological degree of the unordered configuration spaces of a graph; as long as there is at least one essential vertex, this maximal degree is realized for configuration spaces of sufficiently many points (uniformly, twice the number of essential vertices is always enough to realize the maximal degree).
Unstable values are indicated in bold.
All values not explicitly included in the tables of data are stable and are calculated by the indicated stable polynomial value.

\subsection*{Notes about the computation}
The core computation of a presentation for the homology and all Poincar\'e series for $723$ graphs (the isolated vertex has irregularities because it lacks edges and therefore was computed by hand) takes several hours on a 2019 MacBook Pro. The median calculation time is between one and two seconds, while the slowest calculation for a graph in this collection sometimes takes an hour.

Example code for the core computation looks as follows. 
It presents a local two-stage chain complex model for each vertex and then tensors them together, both over the full polynomial ring on the edges.
\begin{adjustwidth}{1cm}{1cm}
\begin{alltt}
-- macaulay script for H_*(B_*(K4))
R = ZZ[e_0, e_1, e_2, e_3, e_4, e_5]
C0 = chainComplex { matrix {{e_4 - e_0, e_5 - e_0}} }
C1 = chainComplex { matrix {{e_1 - e_0, e_2 - e_0}} }
C2 = chainComplex { matrix {{e_3 - e_1, e_5 - e_1}} }
C3 = chainComplex { matrix {{e_3 - e_2, e_4 - e_2}} }
C = C0 ** C1 ** C2 ** C3
H = HH (C)
p0 = hilbertSeries (H_0, Reduce => true)
p1 = hilbertSeries (H_1, Reduce => true)
p2 = hilbertSeries (H_2, Reduce => true)
p3 = hilbertSeries (H_3, Reduce => true)
p4 = hilbertSeries (H_4, Reduce => true)
\end{alltt}
\end{adjustwidth}
\clearpage
\section{Data for small graphs}
\setcounter{subsection}{-1}
\subsection{Data for graphs with 0 essential vertices}
\ 
\vspace{10pt}
\hrule\hrule
\vspace{20pt}
\begin{absolutelynopagebreak}
\begin{center}

\\

Note: values calculated by hand
\\
Betti numbers $\beta_i(B_k(\Gamma))$:
\\
\begin{tabular}{Sl|ScScScSc}
$i\backslash k$ & $0$ & $1$ & $2$ & Poincar\'e series\\\hline
$0$ & $\mathbf{1}$ & $\mathbf{1}$ & $0$
 & $1+t$ \\
\end{tabular}
\begin{tabular}{Sl|Sc}
$i$ & stable polynomial value
\\\hline
$0$ & $0$ \\
\end{tabular}%
\end{center}
\end{absolutelynopagebreak}\newpage
\clearpage
\subsection{Data for graphs with 1 essential vertex}
\ 
\vspace{10pt}
\hrule\hrule
\vspace{20pt}
\begin{absolutelynopagebreak}
\begin{center}
%
\end{center}
\end{absolutelynopagebreak}\newpage
\clearpage
\subsection{Data for graphs with 2 essential vertices}
\ 
\vspace{10pt}
\hrule\hrule
\vspace{20pt}
\begin{absolutelynopagebreak}
\begin{center}
%
\end{center}
\end{absolutelynopagebreak}\newpage
\clearpage
\subsection{Data for graphs with 3 essential vertices}
\ 
\vspace{10pt}
\hrule\hrule
\vspace{20pt}
\begin{absolutelynopagebreak}
\begin{center}
%
\end{center}
\end{absolutelynopagebreak}\newpage
\clearpage
\subsection{Data for graphs with 4 essential vertices}
\ 
\vspace{10pt}
\hrule\hrule
\vspace{20pt}
\begin{absolutelynopagebreak}
\begin{center}
%
\end{center}
\end{absolutelynopagebreak}\newpage
\clearpage
\subsection{Data for graphs with 5 essential vertices}
\ 
\vspace{10pt}
\hrule\hrule
\vspace{20pt}
\begin{absolutelynopagebreak}
\begin{center}
%
\end{center}
\end{absolutelynopagebreak}\newpage
\clearpage
\subsection{Data for graphs with 6 essential vertices}
\ 
\vspace{10pt}
\hrule\hrule
\vspace{20pt}
\begin{absolutelynopagebreak}
\begin{center}
%
\end{center}
\end{absolutelynopagebreak}\newpage
\bibliographystyle{amsalpha}
\begingroup
\raggedright
\bibliography{references}
\endgroup
\end{document}